\documentclass{article}
\usepackage{amsmath}
\usepackage{amsfonts}

\newtheorem{theorem}{Theorem}

\newtheorem{example}[theorem]{Example}

\input{tcilatex}
\begin{document}

\title{Curves with normal planes at constant distance from a fixed point}
\author{Yasemin ALAGOZ \\
{\footnotesize Department of Mathematics, Faculty of Arts and Sciences,
Yildiz Technical University Istanbul-TURKEY}\\
\\
{\footnotesize Email: ygulluk@yildiz.edu.tr, yasemingulluk@yahoo.com}\\
{\footnotesize \ \ \ \ \ \ \ \ \ \ \ \ \ \ }}
\maketitle

\begin{abstract}
We consider the curves whose all normal planes are at the same distance from
a fixed point and obtain some characterizations of them in the 3-dimensional
Euclidean space.
\end{abstract}

\begin{quotation}
\bigskip .

\end{quotation}

\section{Introduction}

Beng-Yen Chen investigated the curves whose position vectors lie in their
rectifying planes in the 3-dimensional Euclidean space, \cite{1}. These
curves are called the rectifying curves which are not plane curves. They can
be equivalently defined as the twisted curves whose osculating planes are at
the same distance from a fixed point. Here we consider the curves whose all
normal planes are at the same distance from a fixed point and obtain some
characterizations of them in the 3-dimensional Euclidean space.

We first review the results on the rectifying curves in the 3-dimensional
Euclidean space, \cite{1}.

Let $C$ be a unit-speed curve with the position vector $\mathbf{r(}s\mathbf{)%
}$. Therefore $s$ is the natural parameter of $C$. Frenet formulas follow as

\begin{equation}
\mathbf{t}^{\prime }=\kappa \mathbf{n},\ \ \ \mathbf{n}^{\prime }=-\kappa 
\mathbf{t}+\tau \mathbf{b},\ \ \ \mathbf{b}^{\prime }=-\tau \mathbf{n}
\label{11}
\end{equation}%
\newline
where$\ \mathbf{t,n}$ and $\mathbf{b}$ are Frenet vectors and $\kappa $ the
curvature and $\tau $ the torsion of the curve.\newline
A rectifying curve can be defined by the equation

\begin{equation}
<\mathbf{r,n>}=0  \label{12}
\end{equation}%
\newline
By means of (\ref{11}) we can write it as

\begin{equation}
<\mathbf{r,t}^{\prime }>=0
\end{equation}%
Accordingly $<\mathbf{r,t>}^{\prime }\mathbf{=}<\mathbf{r,t}^{\prime }%
\mathbf{>+}<\mathbf{t,t>=}1$ and so 
\begin{equation}
<\mathbf{r,t>=}\int_{0}^{s}1ds=s  \label{14}
\end{equation}%
\newline
Integrating \ (\ref{14}) we have

\begin{equation}
<\mathbf{r,r>=}\rho ^{2}=s^{2}+c  \label{15}
\end{equation}%
\newline
where $c\ $is a constant. Without loss the generality we can write (\ref{15}%
) as

\begin{equation}
\rho ^{2}=s^{2}+a^{2}  \label{16}
\end{equation}%
\newline
where $a$ is a positive number.\newline
(\ref{12}) , (\ref{14}) and (\ref{16}) imply that the position vector can be
written as 
\begin{equation}
\mathbf{r}=s\mathbf{t}+a\mathbf{b}  \label{17}
\end{equation}%
\newline
Differentiating (\ref{17}) we have%
\begin{equation*}
\mathbf{t=t+}s\mathbf{(}\kappa \mathbf{n)+}a\mathbf{(-}\tau \mathbf{n)}
\end{equation*}%
and so%
\begin{equation}
s\kappa =a\tau  \label{18}
\end{equation}%
This implies that a rectifying curve is a twisted curve, that is $\tau \neq
0.$\newline
On the other hand, from (\ref{17}) we have%
\begin{equation}
<\mathbf{r,b>}=a  \label{19}
\end{equation}%
\newline
So if a curve is a rectifying curve, then its osculating planes are at the
same distance from the origin. Conversely, let us assume that osculating
planes of a twisted curve are at the same distance from the origin.
Differentiating (\ref{19}) we get (\ref{12}) and so it is a rectifying
curve. Therefore we have the following property:

\textit{A twisted curve is a rectifying curve if and only if its osculating
planes are at the same distance from a fixed point.}

On the other hand, if a unit spherical curve $C_{1}$ is defined by the
equation

\begin{equation}
\mathbf{r}_{1}\mathbf{(}s\mathbf{)}=\frac{\mathbf{r(}s\mathbf{)}}{\sqrt{%
s^{2}+a^{2}}}  \label{20}
\end{equation}%
\newline
where $\mathbf{r(}s\mathbf{)}$ is the position vector of a rectifying curve $%
C$ with the natural parameter $s$, it can be show that the position vector
of a rectifying curve $C$ can be represented as\newline

\begin{equation}
\mathbf{r(}s_{1}\mathbf{)}=(a\sec s_{1})\mathbf{r}_{1}\mathbf{(}s_{1}\mathbf{%
)}  \label{21}
\end{equation}%
\newline
where $s_{1}\ $is the natural parameter of the spherical curve $C_{1}.$ Then
the tangent unit vector of the curve $C$ can be written as

\begin{equation}
\mathbf{t(}s_{1}\mathbf{)}=(\sin s_{1})\mathbf{r}_{1}\mathbf{(}s_{1}\mathbf{%
)+(\cos }s_{1}\mathbf{)\mathbf{t}}_{1}\mathbf{\mathbf{(}}s_{1}\mathbf{%
\mathbf{)}}  \label{22}
\end{equation}%
\newline
where $\mathbf{\mathbf{t}}_{1}\mathbf{\mathbf{(}}s_{1}\mathbf{\mathbf{)}}$
is the tangent unit vector of the spherical curve $C_{1}.\ $Differentiating
the last equation we have

\begin{equation}
\kappa \mathbf{n}a\sec ^{2}s_{1}=\mathbf{(\cos }s_{1}\mathbf{)}(\mathbf{r}%
_{1}+\kappa _{1}\mathbf{n}_{1})  \label{23}
\end{equation}%
\newline
$\kappa _{1}$ is the curvature of the spherical curve $C_{1}.$\newline
From (\ref{23}) we find the following relation between the curvatures of the
curves $C\ $and $C_{1}$ associated with each other: 
\begin{equation}
a^{2}\kappa ^{2}=\mathbf{(\cos }s_{1}\mathbf{)}^{6}(\kappa _{1}^{2}-1)
\label{24}
\end{equation}%
\newline
The last relation can be written, in terms of the natural parameter $s$ of
the rectifying curve, as

\begin{equation}
(s^{2}+a^{2})^{3}\kappa ^{2}=a^{4}(\kappa _{1}^{2}-1)  \label{25}
\end{equation}

\section{Curves with normal planes at the same distance from a fixed point}

Let us assume that all normal planes of a curve $C$ with the position vector$%
\ \mathbf{r=r(}s\mathbf{)}$\textbf{\ }are at the same distance from a fixed
point where $s$ is the natural parameter of $C$. We denote Frenet vectors
and the curvature and the torsion of $C$ by $\mathbf{t,n}$ and $\mathbf{b}$
and $\kappa $ and $\tau $ respectively. We can choose the fixed point as the
origin. Therefore our condition becomes

\begin{equation}
<\mathbf{r,t>}=c_{1}=\text{constant}  \label{31}
\end{equation}%
\newline
By the integration we find that

\begin{equation}
<\mathbf{r,r>=}\rho _{1}^{2}=2c_{1}s+c_{2}  \label{32}
\end{equation}%
\newline
where $c_{2}$ is a constant.\newline
Since the case of $c_{1}=0$ corresponds to a spherical curve, we can assume
that $c_{1}\neq 0$. Then, without loss the generality, we can write (\ref{32}%
) as

\begin{equation}
<\mathbf{r,r>=}\rho _{1}^{2}=4\epsilon c^{2}s,\ \ \ \ \ \epsilon =\text{%
sign(s)}  \label{33}
\end{equation}%
\newline
where $c$ is a positive constant. Then (\ref{31}) reduces to

\begin{equation}
<\mathbf{r,t>}=2\epsilon c^{2}  \label{34}
\end{equation}%
\newline
(\ref{33}) and (\ref{34}) imply 
\begin{equation}
\epsilon s\geq c^{2}  \label{35}
\end{equation}%
Differentiating (\ref{34}) we have

\begin{equation}
\kappa <\mathbf{r,n}>=-1  \label{36}
\end{equation}%
\newline
or

\begin{equation}
<\mathbf{r,n>}=-R  \label{37}
\end{equation}%
\newline
where $R=1/\kappa $ is the radius of curvature.\newline
From (\ref{34}) and (\ref{37}) we can write the position vector as 
\begin{equation}
\mathbf{r=}2\epsilon c^{2}\mathbf{t}-R\mathbf{n}+A\mathbf{b}  \label{38}
\end{equation}%
Differentiating (\ref{38}) we have

\begin{equation}
\mathbf{t}=\mathbf{t}+(2\epsilon c^{2}\kappa -R^{\prime }-A\tau )\mathbf{n}%
+(A^{\prime }-R\tau )\mathbf{b}  \label{39}
\end{equation}%
\newline
and

\begin{equation}
A\tau =2\epsilon c^{2}\kappa -R^{\prime },\ \ \ \ \ \ \ \ A^{\prime }=R\tau
\label{40}
\end{equation}%
\newline
Therefore we first have%
\begin{equation}
RR^{\prime }=2\epsilon c^{2}\Longleftrightarrow \tau =0
\end{equation}%
\newline
So the curve is a plane curve if and only if 
\begin{equation}
RR^{\prime }=2\epsilon c^{2}
\end{equation}%
\newline
Without the generality, for case of a plane curve we can write%
\begin{equation}
\mathbf{r=}2\epsilon c^{2}\mathbf{t+(}-R)\mathbf{n}  \label{401}
\end{equation}%
\newline
Then from (\ref{33}) and (\ref{401}) 
\begin{equation}
R^{2}=4c^{2}(\epsilon s-c^{2})
\end{equation}%
\newline
Accordingly the natural equations of the curve are

\begin{equation}
\kappa =\frac{1}{2c\sqrt{\epsilon s-c^{2}}},\ \ \ \ \ \ \ \tau =0
\label{404}
\end{equation}%
We can now assume that our curve is a twisted curve, that is $\tau \neq 0$.
Then%
\begin{equation}
A=2\epsilon c^{2}\kappa T-R^{\prime }T,\ \ \ \ \ \ \ \ A^{\prime }=R\tau
\label{402}
\end{equation}%
where $T=\frac{1}{\tau }\ $is the radius of torsion. Hence (\ref{38}) can be
written as%
\begin{equation}
\mathbf{r=(}2\epsilon c^{2})\mathbf{t+(}-R)\mathbf{n}+(2\epsilon c^{2}\kappa
T-R^{\prime }T)\mathbf{b}  \label{403}
\end{equation}%
Therefore, according to (\ref{33}) we have%
\begin{equation}
4\epsilon c^{2}s=4c^{4}+R^{2}+(2\epsilon c^{2}\kappa T-R^{\prime }T)^{2}
\label{311}
\end{equation}%
On the other hand, from (\ref{402}) we can write\newline
\begin{equation}
R\tau +(R^{\prime }T)^{\prime }-2\epsilon c^{2}(\kappa T)^{\prime }=0
\label{312}
\end{equation}%
\newline
We can show that (\ref{311}) and (\ref{312}) are equivalent equations. In
fact, the equation (\ref{312}) can be written as

\begin{equation}
\lbrack R\tau +(R^{\prime }T)^{\prime }-2\epsilon c^{2}(\kappa T)^{\prime
}][RR^{\prime }-2\epsilon c^{2}]=0  \label{313}
\end{equation}%
\newline
because of $RR^{\prime }-2\epsilon c^{2}=0$, corresponds to the plane curve.
Then we have 
\begin{equation}
2\epsilon c^{2}=[R\tau +(R^{\prime }T)^{\prime }-2\epsilon c^{2}(\kappa
T)^{\prime }]R^{\prime }T+4c^{4}\kappa T(\kappa T)^{\prime }-2\epsilon
c^{2}\kappa T(R^{\prime }T)^{\prime }
\end{equation}%
and%
\begin{equation}
4c^{2}=2RR^{\prime }+2(2\epsilon c^{2}\kappa T-R^{\prime }T)[2\epsilon
c^{2}(\kappa T)^{\prime }-(R^{\prime }T)^{\prime }]
\end{equation}%
\newline
This implies that%
\begin{equation}
4\epsilon c^{2}=(R^{2})^{\prime }+[(2\epsilon c^{2}\kappa T-R^{\prime
}T)^{2}]^{\prime }  \label{314}
\end{equation}%
\newline
Therefore we have%
\begin{equation*}
R^{2}+(2\epsilon c^{2}\kappa T-R^{\prime }T)^{2}=\int_{c^{2}}^{\epsilon
s}4c^{2}dt=4\epsilon c^{2}s-4c^{4}
\end{equation*}%
It is obvious that using the equation (\ref{311}) we obtain the equation (%
\ref{312}).

\section{Spherical curves associated with a curve with normal planes at
constant distance from a point}

Let us define a unit spherical curve $C_{1}$ by the equation

\begin{equation}
\mathbf{r}_{1}\mathbf{(}s\mathbf{)}=\frac{\mathbf{r(}s\mathbf{)}}{2c\sqrt{%
\epsilon s}}  \label{41}
\end{equation}%
where $\mathbf{r}(s)$ is the position vector of a curve with normal planes
at constant distance from a fixed point.\newline
So we have

\begin{equation}
\mathbf{r(}s\mathbf{)}=2c\sqrt{\epsilon s}\mathbf{r}_{1}\mathbf{(}s\mathbf{)}
\label{42}
\end{equation}%
\newline
Differentiating we get

\begin{equation}
\mathbf{t}=\frac{\epsilon c}{\sqrt{\epsilon s}}\mathbf{r}_{1}\mathbf{(}s%
\mathbf{)}+(2c\sqrt{\epsilon s})\mathbf{r}_{1}^{\prime }\mathbf{(}s\mathbf{)}
\end{equation}%
\newline
From the last equation we find

\begin{equation}
|\mathbf{r}_{1}^{\prime }\mathbf{(}s\mathbf{)}|=\frac{\sqrt{\epsilon s-c^{2}}%
}{2c\epsilon s}=\frac{ds_{1}}{ds}  \label{43}
\end{equation}%
\newline
\newline
for the speed of the spherical curve. And so the natural parameter of the
unit spherical curve$\ C_{1}$, from $s_{1}=\epsilon
\dint\limits_{c^{2}}^{\epsilon s}\frac{\sqrt{t-c^{2}}}{2ct}dt$ is obtained as

\begin{equation}
s_{1}=\epsilon (\frac{\sqrt{\epsilon s-c^{2}}}{c}-\arctan \frac{\sqrt{%
\epsilon s-c^{2}}}{c})  \label{44}
\end{equation}%
\newline
\newline
Since

\begin{equation}
\frac{ds_{1}}{ds}=\frac{\sqrt{\epsilon s-c^{2}}}{2c\epsilon s}>0  \label{45}
\end{equation}%
\newline
\newline
there exists a function $s(s_{1})$ which satisfies the equation (\ref{44})
Therefore (\ref{41}) can be written as

\begin{equation}
\mathbf{r}_{1}(s(s_{1}))=\frac{\mathbf{r}(s(s_{1}))}{2c\sqrt{\epsilon
s(s_{1})}}  \label{46}
\end{equation}%
\newline
So for a given curve $C$ with normal planes at constant distance from a
fixed point, whose position vector is$\ \mathbf{r}(s_{1})=\mathbf{r}%
(s(s_{1}))$, we have a unit spherical curve $C_{1}$ whose position vector $%
\mathbf{r}_{1}(s_{1})=\mathbf{r}_{1}(s(s_{1}))$ is defined by (\ref{46})
with the natural parameter $s_{1}$. We call $C_{1}\ $the unit spherical
curve associated with the curve $C$ with normal planes at constant distance
from a fixed point.

Now let us consider a curve $C$ defined by

\begin{equation}
\mathbf{r}(s_{1})=(2c\sqrt{\epsilon s(s_{1})}\mathbf{r}_{1}(s_{1})
\label{47}
\end{equation}%
\newline
where$\ \mathbf{r}_{1}(s_{1})$ is the position vector of a unit spherical
curve $C_{1}\ $with the natural parameter $s_{1}$ and $s(s_{1})$ is defined
by (\ref{44}) and $c$ is constant. Let $\mathbf{t}_{1}\mathbf{,\ n}_{1}%
\mathbf{,\ b}_{1}$ and $\kappa _{1}$ and$\ \tau _{1}$ be Frenet vectors and
the curvature and the torsion of $C_{1}$ respectively.\newline
According to (\ref{44})

\begin{equation}
s^{\prime }=\frac{ds}{ds_{1}}=\frac{2c\epsilon s}{\sqrt{\epsilon s-c^{2}}}
\label{48}
\end{equation}%
\newline
Differentiating (\ref{47}) with respect to $s_{1}$ we have

\begin{equation}
\mathbf{r}^{\prime }=\frac{d\mathbf{r}}{ds_{1}}=2\epsilon c^{2}\frac{\sqrt{%
\epsilon s}}{\sqrt{\epsilon s-c^{2}}}\mathbf{r}_{1}+(2c\sqrt{\epsilon s})%
\mathbf{t}_{1}  \label{49}
\end{equation}%
\newline
Since $<\mathbf{r}_{1},\mathbf{r}_{1}>=1,\ <\mathbf{r}_{1},\mathbf{t}_{1}>=0$
and $<\mathbf{t}_{1},\mathbf{t}_{1}>=1$ from (\ref{49}) we have

\begin{equation*}
\left\vert \mathbf{r}^{\prime }\right\vert =\frac{2c\epsilon s}{\sqrt{%
\epsilon s-c^{2}}}
\end{equation*}%
\newline
So the unit tangent vector of $C$ is found as

\begin{equation}
\mathbf{t(}s_{1}\mathbf{)}=\frac{\epsilon c}{\sqrt{\epsilon s}}\mathbf{r}%
_{1}(s_{1})+\frac{\sqrt{\epsilon s-c^{2}}}{\sqrt{\epsilon s}}\mathbf{t}%
_{1}(s_{1})  \label{410}
\end{equation}%
\newline
Since $<\mathbf{r},\mathbf{r}_{1}>=2c\sqrt{\epsilon s}\ $and $<\mathbf{r,t}%
_{1}>=0$

\begin{equation}
<\mathbf{r,t>}=2\epsilon c^{2}  \label{411}
\end{equation}%
\newline
This means that the curve $C$ is a curve with normal planes at constant
distance from a fixed point. We call $C\ $the curve with normal planes at
constant distance from a fixed point associated with the unit spherical
curve $C_{1}$.\newline
Let us now differentiate (\ref{410}) with respect to$\ s$. We have

\begin{equation}
\kappa \mathbf{n}=-\frac{c}{2\epsilon s\sqrt{\epsilon s}}\mathbf{r}_{1}+%
\frac{1}{2\sqrt{\epsilon s}\sqrt{\epsilon s-c^{2}}}\mathbf{t}_{1}+\frac{%
\epsilon s-c^{2}}{2\epsilon cs\sqrt{\epsilon s}}\kappa _{1}\mathbf{n}_{1}
\label{412}
\end{equation}%
\newline
Using\textbf{\ }$<\mathbf{n,n>}=1,<\mathbf{r}_{1},\mathbf{r}_{1}>=1,<\mathbf{%
t}_{1},\mathbf{t}_{1}>=1,<\mathbf{n}_{1},\mathbf{n}_{1}>=1,<\mathbf{r}_{1},%
\mathbf{t}_{1}>=0,$ $<\mathbf{t}_{1},\mathbf{n}_{1}>=0$ and $\kappa _{1}<%
\mathbf{r}_{1},\mathbf{n}_{1}>=-1$, from (\ref{412}) we obtain the following
relation between curvatures of the curves $C$ and $C_{1}$ associated with
each other:

\begin{equation}
\kappa ^{2}=\frac{3s^{2}-3c^{2}\epsilon s+c^{4}}{4\epsilon s^{3}(s-c^{2})}+%
\frac{(\epsilon s-c^{2})^{2}}{4c^{2}\epsilon s^{3}}\kappa _{1}^{2}
\label{413}
\end{equation}%
\newline
or%
\begin{equation}
\kappa ^{2}=\frac{1}{4c^{2}(\epsilon s-c^{2})}+\frac{(\epsilon s-c^{2})^{2}}{%
4c^{2}\epsilon s^{3}}(\kappa _{1}^{2}-1)  \label{4131}
\end{equation}%
\newline
Since $\kappa _{1}<\mathbf{r}_{1},\mathbf{n}_{1}>=-1,$the curvature of the
unit spherical curve is not smaller than $1,\ $that is $\kappa _{1}\geq 1\ .$

\begin{example}
\end{example}

The equation (\ref{4131}) implies that the curvature of the unit spherical
curve associated with the plane curve given by the natural equations (\ref%
{404}) is $\kappa _{1}=1.$ This means that the spherical curve is a great
circle of the unit sphere. Hence we can obtain the cartesian equations of
the plane curve associated with a great circle of the unit sphere using the
equation (\ref{42}). In fact, we can choose the equation of a great circle
of the unit sphere as 
\begin{equation}
\mathbf{r}_{1}=(\cos s_{1},\sin s_{1},0)
\end{equation}%
\newline
Since $s_{1}=\epsilon (\frac{\sqrt{\epsilon s-c^{2}}}{c}-\arctan \frac{\sqrt{%
\epsilon s-c^{2}}}{c}),$%
\begin{equation*}
\cos s_{1}=\cos \frac{\sqrt{\epsilon s-c^{2}}}{c}\frac{c}{\sqrt{\epsilon s}}%
+\epsilon \sin \frac{\sqrt{\epsilon s-c^{2}}}{c}\frac{\sqrt{\epsilon s-c^{2}}%
}{\sqrt{\epsilon s}}
\end{equation*}%
and%
\begin{equation*}
\sin s_{1}=\epsilon \sin \frac{\sqrt{\epsilon s-c^{2}}}{c}\frac{c}{\sqrt{%
\epsilon s}}-\cos \frac{\sqrt{\epsilon s-c^{2}}}{c}\frac{\sqrt{\epsilon
s-c^{2}}}{\sqrt{\epsilon s}}
\end{equation*}%
Then from (\ref{42}) we obtain the Cartesian equations of the plane curve
with normal planes at constant distance from the origin as%
\begin{eqnarray}
x &=&2c^{2}\cos \frac{\sqrt{\epsilon s-c^{2}}}{c}+2\epsilon c\sqrt{\epsilon
s-c^{2}}\sin \frac{\sqrt{\epsilon s-c^{2}}}{c}  \notag \\
y &=&2\epsilon c^{2}\sin \frac{\sqrt{\epsilon s-c^{2}}}{c}-2c\sqrt{\epsilon
s-c^{2}}\cos \frac{\sqrt{\epsilon s-c^{2}}}{c}  \label{385} \\
z &=&0  \notag
\end{eqnarray}%
This curve is a plane curve with normal planes at a distance $2c^{2}$ from
the origin. This means that all normal lines of the curve are at a distance $%
2c^{2}$ from the origin. So it is an involute of the circle of radius $%
2c^{2} $ centered at the origin. Since it is an involute of a plane curve,
it is also a plane curve \cite[p,88]{2}. In fact, the position vector of an
involute $C$ of the curve $C_{2}$ whose position vector is$\ \mathbf{r}%
_{2}(s_{2})$can be written as%
\begin{equation}
\mathbf{r=r}_{2}(s_{2})+(c_{2}-s_{2})\mathbf{t}_{2}(s_{2})  \label{3850}
\end{equation}%
where $c_{2}\ $is constant and $\mathbf{t}_{2}(s_{2})$ is unit tangent
vector to the curve $C_{2}$, \cite[p.69]{2}, \cite[p.40]{4}, \cite[p.99]{5}.
Then%
\begin{equation}
\mathbf{r}^{\prime }\mathbf{=}(c_{2}-s_{2})\kappa _{2}\mathbf{n}_{2}
\end{equation}%
where $\kappa _{2}$ is the curvature of $C_{2}$ and $\mathbf{n}_{2}$ the
principal normal vector. Therefore the unit tangent vector t to $C$ can be
written as 
\begin{equation}
\mathbf{t=-}\epsilon _{2}\mathbf{n}_{2},\ \ \ \ \ \ \mathbf{-}\epsilon _{2}=%
\text{sign}(c_{2}-s_{2})  \label{3851}
\end{equation}%
Since $C_{2}$ is a circle of radius $2c^{2}$ centered at the origin,%
\begin{equation}
\mathbf{r}_{2}=(2c^{2}\cos \frac{s_{2}}{2c^{2}},2c^{2}\sin \frac{s_{2}}{%
2c^{2}},0)
\end{equation}%
\begin{equation}
\mathbf{t}_{2}=(-\sin \frac{s_{2}}{2c^{2}},\cos \frac{s_{2}}{2c^{2}},0)
\end{equation}%
\begin{equation}
\mathbf{n}_{2}=(-\cos \frac{s_{2}}{2c^{2}},-\sin \frac{s_{2}}{2c^{2}},0)
\label{3852}
\end{equation}%
So$\ <\mathbf{r}_{2},\mathbf{t}_{2}>=0$, $<\mathbf{r}_{2},\mathbf{n}%
_{2}>=-2c^{2}$. Then from (\ref{3850}) and (\ref{3851}) we have 
\begin{equation}
<\mathbf{r},\mathbf{t}>=2\epsilon _{2}c^{2}
\end{equation}%
This means that the involute $C$ of the circle $C_{2}$ is a plane curve with
normal planes at a distance $2c^{2}$ from the origin.

\FRAME{dtbpFU}{2.2796in}{2.2796in}{0pt}{\Qcb{Figure 1}}{}{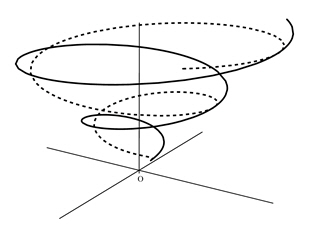}{%
\special{language "Scientific Word";type "GRAPHIC";maintain-aspect-ratio
TRUE;display "USEDEF";valid_file "F";width 2.2796in;height 2.2796in;depth
0pt;original-width 2.2399in;original-height 2.2399in;cropleft "0";croptop
"1";cropright "1";cropbottom "0";filename
'../040216/untitled2.png';file-properties "XNPEU";}}

Figure 1 illustrates the plane curve (\ref{385}) with normal planes at a
distance $2c^{2}=2(1)^{2}=2$ from the origin (the solid line for $\epsilon
=1 $, the dashed line for $\epsilon =-1$). It is an involute of the circle
of radius $2c^{2}=2$ centered at the origin.\newline
Let us note the equations (\ref{385}) can be also obtained using the
equations (\ref{3850})-(\ref{3852}) or using (\ref{385}) from the equations
of a plane curve given by%
\begin{equation}
\mathbf{r}=\mathbf{r}(s)=(x(s),y(s),0)=(\dint \cos (\dint \kappa
(s)ds)ds,\dint \sin (\dint \kappa (s)ds)ds,0)
\end{equation}%
where $s$ is the natural parameter of the curve, \cite[p.87]{2}, \cite[p.28]%
{3}, \cite[p.99]{5}. \newline
In the following we give an example of a twisted curve with normal planes at
constant distance from a fixed point.\newline

\begin{example}
\ 
\end{example}

Let us choose the unit spherical curve $C_{1}$ as the circle of radius $%
\frac{1}{\sqrt{2}}$ given by the equation

\begin{equation}
\mathbf{r}_{1}=(\frac{1}{\sqrt{2}}\cos (\sqrt{2}s_{1})\mathbf{,}\frac{1}{%
\sqrt{2}}\sin (\sqrt{2}s_{1}),\frac{1}{\sqrt{2}})  \label{414}
\end{equation}%
\newline
Then the curve $C$ associated with $C_{1}$ is given by the equation%
\begin{equation}
\mathbf{r(}s\mathbf{)}=(c\sqrt{2}\sqrt{\epsilon s}\cos (\sqrt{2}s_{1})%
\mathbf{,}c\sqrt{2}\sqrt{\epsilon s}\sin (\sqrt{2}s_{1}),c\sqrt{2}\sqrt{%
\epsilon s})  \label{4141}
\end{equation}%
where $s_{1}=\epsilon (\frac{\sqrt{\epsilon s-c^{2}}}{c}-\arctan \frac{\sqrt{%
\epsilon s-c^{2}}}{c}).$\newline
Since the curvature of the circle $C_{1}$ is $\kappa _{1}=\sqrt{2},\ $(\ref%
{4131}) implies that%
\begin{equation}
\kappa ^{2}=\frac{1}{4c^{2}(\epsilon s-c^{2})}+\frac{(\epsilon s-c^{2})^{2}}{%
4c^{2}\epsilon s^{3}}=\frac{\epsilon s^{3}+(\epsilon s-c^{2})^{3}}{%
4c^{2}\epsilon s^{3}(\epsilon s-c^{2})}
\end{equation}%
\newline
Figure $2$ illustrates the twisted curve (\ref{4141}) with normal planes at
a distance $2c^{2}=2(1)^{2}=2$ from the origin (the solid line for $\epsilon
=1$, the dashed line for $\epsilon =-1$).\FRAME{dtbpFU}{2.1188in}{1.6094in}{%
0pt}{\Qcb{Figure 2}}{}{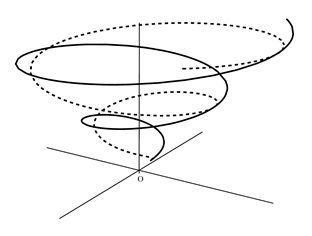}{\special{language "Scientific
Word";type "GRAPHIC";maintain-aspect-ratio TRUE;display "USEDEF";valid_file
"F";width 2.1188in;height 1.6094in;depth 0pt;original-width
2.0799in;original-height 1.5731in;cropleft "0";croptop "1";cropright
"1";cropbottom "0";filename '../040216/untitled4.png';file-properties
"XNPEU";}}\newline

\section{Curves of constant curvature with normal planes at the same
distance from a fixed point}

If $\kappa =\frac{1}{r}=$constant, (\ref{403}) reduce to%
\begin{equation}
\mathbf{r}=(2\epsilon c^{2})\mathbf{t}+(-r)\mathbf{n}+2\frac{\epsilon c^{2}}{%
r}T\mathbf{b}  \label{4151}
\end{equation}%
From (\ref{4151}) we have $<\mathbf{r,n>}=-r=$constant. This means that
rectifying planes of a curve of constant curvature with normal planes at
constant distance from a fixed point are also at constant distance from the
same point.

Conversely if normal planes and rectifying planes of a curve are at constant
distance from a fixed point, it is a curve of constant curvature. In fact by
above conditions we can write

\begin{equation}
<\mathbf{r,t>}=c_{1}=\text{constant}  \label{415}
\end{equation}%
\newline
and

\begin{equation}
<\mathbf{r,n>}=c_{2}=\text{constant}  \label{416}
\end{equation}%
\newline
Differentiating (\ref{415}) we have

\begin{equation*}
1+<\mathbf{r},\kappa \mathbf{n>}=0
\end{equation*}%
\newline
and so

\begin{equation}
<\mathbf{r,n>}=-R  \label{417}
\end{equation}%
\newline
Because of (\ref{416}) and (\ref{417}) the curve is a curve of constant
curvature.

On the other hand for a curve of constant curvature with normal planes at
constant distance from a fixed point (\ref{311}) reduces to

\begin{equation*}
4\epsilon c^{2}s=4c^{4}+r^{2}+\frac{4c^{4}}{r^{2}}T^{2}
\end{equation*}%
\newline
Therefore the natural equations of such a curve are

\begin{equation*}
R=r=\text{constant},\ \ \ \ \ \ T^{2}=r^{2}\frac{2\epsilon hs-r^{2}-h^{2}}{%
h^{2}}
\end{equation*}%
\newline
where $h=2c^{2}$ is the distance of the normal planes from the origin and $r$
is the distance of the rectifying planes from the origin.


\begin{thebibliography}{9}
\bibitem{1} Chen BY. When does the position vector of a space curve always
lie in its rectifying plane?. Am Math Mon 110 2003; VOL: 147-152

\bibitem{3} Eisenhart LP. A Treatise on The Differential Geometry of Curves
and Surfaces. New York: Dover Publications, Inc, 1909.

\bibitem{2} Goetz A. Introduction to Differential Geometry. California,
London, Ontario: Addison Wesley Publication Company, Reading, Massachusetts,
1970.

\bibitem{5} Hsiung CC. A First Course in Differential Geometry. Cambridge,
MA 02238-2872 USA: International Press, P.O. Box 2872 , 1997

\bibitem{4} Struik DJ. Lectures on Classical Differential Geometry.
Cambridge,42, Mass: Addison-Wesley Press, Inc, 1950.
\end{thebibliography}
\end{document}